\documentclass[twoside,11pt]{article}
\title{} \author{} \date{}
\usepackage{latexsym,amssymb,times}%,showkeys}
\input amssym.def
\newtheorem{te}{Theorem}[section]
\newtheorem{prop}[te]{Proposition}
\newtheorem{cor}[te]{Corollary}
\newtheorem{fac}[te]{Fact}
\newtheorem{lem}[te]{Lemma}
\newtheorem{cla}[te]{Claim}

\newtheorem{ex}[te]{Example}

\def\dok{\noindent{\bf Proof. }}
\def\kdok{\hfill $\Box$ \par \vspace*{2mm} }
\def\a{\alpha}
\def\b{\beta}
\def\g{\gamma}
\def\d{\delta}
\def\f{\varphi}
\def\o{\omega}
\def\k{\kappa}

\def\t{{\mathfrak t}}
\def\h{{\mathfrak h}}
\def\A{{\mathbb A}}
\def\B{{\mathbb B}}
\def\S{{\mathbb S}}
\def\P{{\mathbb P}}
\def\Q{{\mathbb Q}}
\def\R{{\mathbb R}}
\def\B{{\mathbb B}}
\def\N{{\mathbb N}}
\def\X{{\mathbb X}}
\def\Y{{\mathbb Y}}
\def\L{{\mathbb L}}
\def\H{{\mathcal H}}
\def\CS{{\mathcal S}}
\def\I{{\mathcal I}}
\def\J{{\mathcal J}}
\def\la{\langle}
\def\ra{\rangle}
\def\lex{<_{{\mathrm{lex}}}}
\def\iso{\stackrel{{\mathrm{iso}}}{\longrightarrow}}
\def\type{\mathop{\mbox{type}}\nolimits}
\def\Lim{\mathop{\mbox{Lim}}\nolimits}
\def\Fin{\mathop{\rm Fin}\nolimits}
\def\sm{\mathop{\rm sm}\nolimits}
\def\sq{\mathop{\rm sq}\nolimits}
\def\rp{\mathop{\rm rp}\nolimits}
\def\supp{\mathop{\rm supp}\nolimits}
\def\red{\mathop{\rm red}\nolimits}
\def\dom{\mathop{\rm dom}\nolimits}
\begin{document}
\thispagestyle{plain}
\begin{center}
           {\large \bf {\uppercase{Forcing with copies of countable ordinals}}}
\end{center}
\begin{center}
{\bf Milo\v s S.\ Kurili\'c}\\[1mm]
  {\small Department of Mathematics and Informatics, University of Novi Sad, \\
         Trg Dositeja Obradovi\'ca 4, 21000 Novi Sad, Serbia\\[-1mm]
                                     e-mail: milos@dmi.uns.ac.rs
   }
\end{center}
\begin{abstract}
\noindent
Let $\a $ be a countable ordinal and $\P (\a )$ the collection of its subsets isomorphic to $\a$.
We show that the separative quotient of the poset $\la \P (\alpha ), \subset \ra$ is isomorphic to a forcing product
of iterated reduced products of Boolean algebras of the form $P(\o ^\g )/\I _{\o ^\g}$, where $\g\in \Lim \cup \{ 1 \}$
and $\I _{\o ^\g}$ is the corresponding ordinal ideal.
Moreover, the poset $\la \P (\a ), \subset \ra$ is forcing equivalent to a two-step iteration of the form
$(P(\o )/\Fin )^+ \ast \pi$, where $[\o ] \Vdash$ ``$\pi$ is an $\o _1$-closed separative pre-order"
and, if $\h=\o _1$, to $(P(\o )/\Fin )^+ $.
Also we analyze the quotients over ordinal ideals $P(\o ^\d )/\I _{\o ^\d }$ and the corresponding cardinal invariants
$\h _{\o ^\d }$ and $\t _{\o ^\d }$.
\end{abstract}
\section{Introduction}\label{S1}
The posets of the form $\la \P (\X), \subset \ra$, where $\X$ is a relational structure and
$\P(\X)$ the set of (the domains of) its isomorphic substructures, were considered in \cite{Kur},
where a classification of the relations on countable sets related to the forcing-related properties of the corresponding posets of copies is
described. So, defining two structures to be equivalent if the corresponding posets of copies produce the same
generic extensions, we obtain a rough classification of structures which, in general, depends on the properties of the model of
set theory in which we work.

For example, under CH all countable linear orders are partitioned in only two classes. Namely, by  \cite{KurTod},
CH implies that for a non-scattered countable linear order $L$ the poset $\la \P (L), \subset \ra$ is forcing equivalent to
the iteration $\S \ast (P(\check{\o })/\Fin )^+$, where $\S$ is the Sacks forcing.
Otherwise, for scattered orders,  by \cite{Kur1} we have
\begin{te} \rm \label{T4121}
For each countable scattered linear order $L$   the separative quotient of the poset
$\la \P (L), \subset \ra$ is $\o _1$-closed and atomless. Under CH, it is forcing equivalent to
the poset $(P(\o )/\Fin )^+$.
\end{te}
The aim of this paper is to get a sharper picture of countable scattered linear orders in this context
and we concentrate our attention on ordinals $\a <\o _1$. So, in Section \ref{S3} we describe the separative quotient
of the poset $\la \P (\a ), \subset \ra$ and, in Section \ref{S5}, factorize it  as a two-step iteration
$(P(\o )/\Fin )^+ \ast \pi$, where $[\o ] \Vdash$ ``$\pi$ is an $\o _1$-closed separative pre-order"
(which implies that the equality $\h =\o _1$ implies that all posets  $\la \P (\a ), \subset \ra$  are forcing equivalent to $(P(\o )/\Fin )^+$ again).
In Section \ref{S4} we factorize the quotients $P(\o ^\g )/\I _{\o ^\g}$, for $\g\in \Lim $, and, in
Section \ref{S6}, consider the quotients over the ordinal ideals $P(\o ^\d )/\I _{\o ^\d }$, $0<\d <\o _1$, and analyze the corresponding cardinal invariants
$\h _{\o ^\d }$ and $\t _{\o ^\d }$.
\section{Preliminaries}\label{S2}
In this section we recall some  definitions and basic facts used in the paper.

If $\X $ is a relational structure, $X$ its domain and $A\subset X$, then $\A$ will denote the corresponding substructure of $\X$.
Let $\P (\X )=\{ A \subset X : \A \cong \X\}$
and $\I _{\X }= \{ A\subset X : \X \not\hookrightarrow \A \}$.
It is easy to check that $\X$ is an {\it indivisible structure} (that is, for each partition $X=A\cup B$ we have $\X \hookrightarrow \A$, or
$\X \hookrightarrow \B$) iff $\I _{\X }$ is an ideal. We will use the following elementary fact.
\begin{fac}\rm\label{T4056}
Let $\X $ and $\Y $ be relational structures and $f: \X \iso \Y $. Then

(a) $A\in \I _\X \Leftrightarrow f[A]\in \I _\Y$, for each $A\subset X$;

(b) $\la P(X)\setminus \I _\X , \subset \ra  \cong \la P(Y)\setminus \I _\Y , \subset \ra$.
\end{fac}
A linear order $L$ is said to be {\it scattered} iff it does not contain a dense suborder or, equivalently,
iff the rational line, $\Q$, does not embed in $L$.
By $\CS$ we denote the class of all countable scattered linear orders.
A linear order $L$ is said to be {\it additively indecomposable}
iff for each decomposition
$L=L_0 + L_1$ we have $L\hookrightarrow L_0$ or $L\hookrightarrow L_1$. The class $\H$ of {\it hereditarily additively
indecomposable} (or {\it ha-indecomposable}) linear orders is the smallest class of order types of
countable linear orders containing the one element order type, {\bf 1}, and containing the $\o$-sum, $\sum _\o L_i$,
and the $\o ^*$-sum, $\sum _{\o ^*} L_i$, for each sequence $\la L_i :i\in \o \ra$ in $\H$ satisfying
\begin{equation}\label{EQ4100}
\forall i\in \o \;\; |\{ j\in \o : L_i \hookrightarrow L_j \}|=\aleph _0 .
\end{equation}
\begin{fac}  \rm \label{T4100}
(Laver, \cite{Lav}) $\H \subset \CS$. If $L\in \CS$, then $L\in \H$ iff $L$ is additively indecomposable (see also \cite{Rosen}, p.\ 196 and p.\ 201).
\end{fac}
\begin{fac}\rm\label{T4408}
(See \cite{Kur1})
(a) Let $L=\sum _\o L_i \in \H$, where $\la L_i :i\in \o \ra$ is a sequence in $\H$ satisfying (\ref{EQ4100}). Then
$A\subset L$ contains a copy of $L$ iff
for each $ i, m\in \o $ there is finite $K \subset \o \setminus m$ such that
$L_i \hookrightarrow \bigcup _{j\in K }L_j \cap A $.

(b) Let $L=\sum _{i\leq n } L_i$, where  $L_i\in \H$ are $\o$-sums
of sequences in $\H$ satisfying (\ref{EQ4100})
and $L_i +L_{i+1}\not\in \H$, for $i<n$. Then
$%$\textstyle
\la \P (L), \subset \ra \cong \prod _{i\leq n } \la \P (L_i), \subset \ra .
$%$
\end{fac}
If $\la A, <\ra$ is a well ordering, $\type \la A, <\ra$  denotes the unique ordinal isomorphic to $\la A, <\ra$.
The product of ordinals $\a$ and $\b$ is the ordinal $\a \b =\type \la \b \times \a ,\lex \ra$, where
$\lex$ is the lexicographic order on the product $\b \times \a$ defined by
$\la \xi , \zeta \ra \lex \la \xi ',\zeta '\ra \Leftrightarrow \xi < \xi ' \lor (\xi = \xi ' \land \zeta < \zeta ')$.
The power $\a ^\b$ is defined recursively by
$\a ^0 =1 , \; \a ^{\b +1}= \a ^\b  \a \mbox{ and } \a ^\g =\sup \{ \a ^\xi : \xi < \g\} , \mbox{ for limit }\g $.
For an ordinal $\a$, instead of $\P (\la \a , \in \ra )$ we will write $\P (\a)$.
\begin{fac}\rm\label{T4400}
For a countable limit ordinal $\a$ the following conditions are equivalent:

(a) $\a$ is indecomposable (i.e.\ $\a $ is not a sum of two smaller ordinals);

(b) $\b + \g < \a$, for each $\b , \g <\a$;

(c) $A\in \P(\a )$ or $\a \setminus A\in \P (\a )$, for each $A\subset \a$;

(d) $\a = \o ^\delta$, for some countable ordinal $\delta >0$;

(e) $\a\in \H$;

(f) $\a$ is an indivisible structure;

(g) $\I _\a =\{ I\subset \a : \a \not\hookrightarrow I\}$ is an ideal in $P(\a )$.
\end{fac}
\dok
For the equivalence of (b), (c) and (d) see \cite{Kun}, p.\ 43.
For (a) $\Leftrightarrow $ (d) see 1.3.6 of \cite{Fra}.
By \cite{Rosen}, p.\ 176, (d) holds iff $\a$ is additively indecomposable which is, by Fact \ref{T4100},
equivalent to (e). (a) $\Leftrightarrow $ (f) is 6.8.1 of \cite{Fra}. (f) $\Leftrightarrow$ (g) is evident.
\hfill $\Box$
\begin{fac}\rm\label{T4401}
For each ordinal $\a$ we have $\P (\a )=P(\a )\setminus \I _\a$. Thus $\P (\o ^\delta )= (\I _{\o ^\delta})^+$.
\end{fac}
\dok
The inclusion ``$\subset$" is trivial. If $\a \hookrightarrow A \subset \a$ then, using the fact that for each
increasing function $f: \a \rightarrow \a$ we have $\b \leq f(\b)$, for each $\b \in \a$, we easily show that
$\type (A)=\a$, which means that $A\in \P (\a)$.
%\hfill $\Box$
\kdok
\noindent
A partial order ${\mathbb P} =\langle  P , \leq \rangle $ is called
{\it separative} iff for each $p,q\in P$ satisfying $p\not\leq q$ there is $r\leq p$ such that $r \perp q$.
The {\it separative modification} of ${\mathbb P}$
is the separative pre-order $\mathop{\rm sm}\nolimits ({\mathbb P} )=\langle  P , \leq ^*\rangle $, where
%\begin{equation}\label{EQ4141}
$p\leq ^* q \Leftrightarrow \forall r\leq p \; \exists s \leq r \; s\leq q $.
%\end{equation}
The {\it separative quotient} of ${\mathbb P}$
is the separative partial order $\mathop{\rm sq}\nolimits  ({\mathbb P} )=\langle P /\!\! =^* , \trianglelefteq \rangle$, where
$p = ^* q \Leftrightarrow p \leq ^* q \land q \leq ^* p\;$ and $\;[p] \trianglelefteq [q] \Leftrightarrow p \leq ^* q $ (see \cite{Jech}).
\begin{fac}  \rm \label{T4042}
Let ${\mathbb P} , {\mathbb Q} $ and ${\mathbb P} _i$, $i\in I$, be partial orderings. Then

(a) ${\mathbb P}$, $\mathop{\rm sm}\nolimits ({\mathbb P})$ and $\mathop{\rm sq}\nolimits  ({\mathbb P})$ are forcing equivalent forcing notions;

(b) ${\mathbb P} \cong {\mathbb Q}$ implies that $\mathop{\rm sm}\nolimits {\mathbb P} \cong  \mathop{\rm sm}\nolimits {\mathbb Q}$ and $\mathop{\rm sq}\nolimits  {\mathbb P} \cong  \mathop{\rm sq}\nolimits  {\mathbb Q}$;

(c) $\mathop{\rm sm}\nolimits (\prod _{i\in I}{\mathbb P} _i) = \prod _{i\in I}\mathop{\rm sm}\nolimits {\mathbb P} _i$ and
 $\mathop{\rm sq}\nolimits  (\prod _{i\in I}{\mathbb P} _i) \cong \prod _{i\in I}\mathop{\rm sq}\nolimits  {\mathbb P} _i$.

\noindent
Let $X$ be an infinite set, ${\mathcal I}\varsubsetneq P(X)$ an ideal and $[X]^{<\o }\subset \I$. Then

(d) $\sm \langle  P(X)\setminus {\mathcal I}, \subset  \rangle = \langle  P(X)\setminus {\mathcal I}, \subset _\I \ra$, where
$A\subset _{{\mathcal I}}B \Leftrightarrow A\setminus B \in {\mathcal I}$.

(e) $\sq \langle  P(X)\setminus {\mathcal I}, \subset  \rangle =\langle  (P(X) / =_{{\mathcal I}})^+ , \leq _{{\mathcal I}} \rangle$, where
$A= _{{\mathcal I}}B \Leftrightarrow A \vartriangle B \in {\mathcal I}$ and
$[A] \leq _{{\mathcal I}} [B] \Leftrightarrow A\setminus B \in {\mathcal I}$. Usually this poset is denoted by $(P(X)/\I )^+$.
\end{fac}
\noindent
Let $\k$ be a regular cardinal. A pre-order $\langle {\mathbb P} , \leq \rangle$ is {\it $\k$-closed} iff for each
$\gamma <\k$ and each sequence $\la p_\a :\a <\gamma\ra$ in $\P$, such that $\a <\beta \Rightarrow p_{\beta}\leq p_\a$,
there is $p\in \P$ such that $p\leq p_\a$, for all $\a <\gamma$.
\begin{fac}\rm\label{T4043}
Let $\k$ be a regular cardinal and $\lambda$ an infinite cardinal. Then

(a) If $\P _i$, $i\in I$, are $\k$-closed pre-orders, then the product $\prod _{i\in I}\P _i$ is $\k$-closed.

(b) If ${\mathfrak c}=\o _1$, then each atomless separative $\o _1$-closed
pre-order of size $\o _1$ is forcing equivalent to $(P(\o )/\Fin)^+$ (and to the collapsing algebra Coll$(\o _1 ,\o _1 )$).

(c) If $\lambda ^{<\k }=\lambda $, then each atomless separative $\k$-closed
pre-order $\P$ of size $\lambda$, such that $1_\P \Vdash |\check{\lambda }|=\check{\k }$,
is forcing equivalent to the collapsing algebra Coll$(\k ,\lambda )$.
\end{fac}
\section{The separative quotient of $\la \P(\a ), \subset \ra $}\label{S3}
\noindent
For a Boolean lattice $\B =\la B, \leq \ra$, by $\rp (\B )$ we will denote the {\it reduced power} $\la B ^\o / \equiv , \leq _\equiv \ra $, where
for $\la b_i \ra , \la c_i \ra \in B^\o$, $\la b_i \ra \equiv \la c_i \ra$ (resp.\ $[\la b_i \ra]_\equiv \leq_ \equiv [\la c_i \ra ]_\equiv $)
iff $b_i=c_i$ (resp.\ $b_i \leq c_i$), for all but finitely many $i\in \o$.
For $n\in \o$ we define the set $\rp ^n(\B )$ by: $\rp ^0(\B) =\B $ and $\rp ^{n+1}(\B )= \rp (\rp ^n (\B ))$.

The aim of this section is to prove the following statement.
\begin{te}\rm\label{T4407}
If $\a =\o ^{\g _n +r_n }s_n + \dots + \o ^{ \g _0 +r_0 }s_0 +k $ is a countable ordinal
presented in the Cantor normal form, where $k\in \o$, $r_i \in \o$, $s_i \in \N$, $\g _i \in \Lim \cup \{ 1 \}$ and
$\g _n +r_n > \dots > \g _0 +r_0$,
then
$$%\textstyle
\sq \la \P (\a ), \subset \ra \cong \prod _{i=0}^n \Big( \Big( \rp ^{r_i}( P(\o ^{\g _i} )/ \I _{\o ^{\g _i} })\Big)^+ \Big)^{s_i} .
$$
\end{te}
A proof of Theorem \ref{T4407} is given at the end of the section.

We remind the reader that, if $\I$ and $\J$ are ideals on the sets $X$ and $Y$ respectively, then their {\it Fubini product}
$\I \times \J$ is the ideal on the set $X\times Y$ defined by
$%$
\I \times \J = \{ A\subset X\times Y :
\{ x\in X : \pi _Y [A \cap (\{ x \} \times Y)] \in \J ^+ \} \in \I \} ,
$%$
where $\pi _Y : X\times Y \rightarrow Y$ is the projection.
In particular, if $X=\o$, $\I =\Fin$ and $L_i =\{ i\} \times Y$, for $i\in \o$, then for $A\subset \o \times Y$ we have
%$$
%A\in \Fin \times Y \Leftrightarrow \pi _Y [A\cap L_i] \in \J ^+ \mbox{ for finitely many } i\in \o .
%$$
\begin{equation}\label{EQ4400}
A\in \Fin \times \J \Leftrightarrow
\exists j\in \o \;\; \forall i\geq j \;\; \pi _Y [A\cap L_i] \in \J   .
\end{equation}
For convenience let us define the sets $\o ^n \times Y$, $n\in \o$, recursively by $\o ^0 \times Y =Y$ and $\o ^{n+1} \times Y =\o \times (\o ^n \times Y)$.
Also we define the ideal $\Fin ^n \times \J$ on the set $\o ^n \times Y$ by:
$\Fin ^0 \times \J =  \J$ and $\Fin ^{n+1} \times \J = \Fin \times ( \Fin ^n \times \J)$. Some parts of the following
lemma are folklore but, for completeness,  we include their proofs.
\begin{lem}\rm\label{T4403}
For each ordinal $1\leq \b <\o _1$ and each $n\in \o$ we have:

(a) $\la \P (\o ^{\b +n }), \subset \ra \cong \la P(\o ^n \times \o ^\b )\setminus (\Fin ^n \times \I _{\o ^\b } ) , \subset  \ra$;

(b) $\I _{\o ^{\b +n }}\cong \Fin ^n \times \I _{\o ^\b }$;

(c) $\sq \la\P (\o ^{\b +n }), \subset \ra \cong (P(\o ^n \times \o ^\b )/(\Fin ^n \times \I _{\o ^\b } ))^+  $;

(d) $P(\o ^n \times \o ^\b )/(\Fin ^n \times \I _{\o ^\b } ) \cong \rp ^n (P(\o ^\b )/ \I _{\o ^\b })$;

(e) $\sq (\P (\o ^{\b +n }), \subset ) \cong (\rp ^n (P(\o ^\b )/ \I _{\o ^\b }))^+ $.
\end{lem}
\dok
For $n=0$ the statement follows from Fact \ref{T4401}. So, in the sequel we prove the statement for $n\in \N$.

Using induction we prove (a) and (b) simultaneously. First we show that
\begin{equation}\label{EQ4411}
\la \P (\o ^{\b +1 }), \subset \ra \cong \la (\Fin  \times \I _{\o ^\b } )^+ , \subset  \ra .
\end{equation}
By the properties  of ordinal multiplication and exponentiation we have
$\la \o ^{\b +1}, \in \ra$ $ = \la \o ^\b \o , \in \ra \cong \la \o \times \o ^\b ,\lex \ra= \L$,
where $\L =\sum _{i\in \o}\L _i$ and, for $i\in \o$, $\L _i =\la L_i , <_i\ra$,
$L_i = \{ i \} \times \o ^\b $ and
$\la i, \xi \ra <_i \la i, \zeta \ra \Leftrightarrow \xi \in \zeta $, for $\xi , \zeta \in \o ^\b$.
So, for the function $f_i : L_i \rightarrow \o ^\b$ defined by $f_i (\la i, \xi \ra) =\xi $ we have
\begin{equation}\label{EQ4403}
f_i = \pi _{\o ^\b }\mid L_i : \la L_i , <_i \ra \iso \la \o ^\b ,\in \ra .
\end{equation}
Since $\la \o ^{\b +1 }, \in \ra \cong \L $,
using Facts \ref{T4401} and \ref{T4056}(b) we obtain $\la \P (\o ^{\b +1}), \subset \ra =\la P (\o ^{\b +1})\setminus \I _{\o ^{\b +1}}, \subset \ra
\cong \la P(\L )\setminus \I _\L , \subset \ra$ so it remains to be shown that
\begin{equation}\label{EQ4405}
\I _\L ^+ = (\Fin \times \I _{\o ^\b })^+ .
\end{equation}
\begin{cla} \rm \label{T4437}
For each $A\subset \o \times \o ^\b $ we have

(i) $A\in \I _\L ^+ \Leftrightarrow \forall j\in \o \;\; \exists K \in [\o \setminus j]^{<\o }\;\; \o ^\b \hookrightarrow \bigcup _{i\in K}L_i \cap A$;

(ii) $A\in (\Fin \times  \I _{\o ^\b }) ^+ \Leftrightarrow \forall j\in \o \;\; \exists i\geq j\;\; \o ^\b \hookrightarrow L_i \cap A$.
\end{cla}
\dok
(i) By (\ref{EQ4403}), for each $i\in\o$ we have $\L _i \cong \o ^\b$ so, by Fact \ref{T4400} we have $\L _i \in \H$ and, clearly, condition (\ref{EQ4100}) is satisfied.
By Fact \ref{T4408}(a),  $A\in \I _\L ^+ $ iff
$\forall j\in \o \;\; \exists K \in [\o \setminus j]^{<\o }\;\; \o ^\b \hookrightarrow \bigcup _{i\in K}L_i \cap A$.

(ii) By (\ref{EQ4400}), $A\not\in \Fin \times \I _{\o ^\b }$ iff for each $j\in \o$ there exists $i\geq j $ such that
$\pi _{\o ^\b }[L_i \cap A]\not\in \I _{\o ^\b}$.
But, by (\ref{EQ4403}) and Fact \ref{T4056}(a) we have: $\pi _{\o ^\b }[L_i \cap A]\not\in \I _{\o ^\b}$ iff $f_i [L_i \cap A]\not\in \I _{\o ^\b}$
iff $L_i \cap A \not\in \I _{\L _i}$ iff $\L _i \hookrightarrow L_i \cap A$ iff $\o ^\b \hookrightarrow L_i \cap A$.
\kdok
By Claim \ref{T4437}, the inclusion ``$\supset$" in (\ref{EQ4405}) is satisfied and we prove ``$\subset$". If $A\in \I _\L ^+$ and $j \in \o$ then,
by Claim \ref{T4437}(i), there are
$K\in [\o \setminus j ]^{<\o }$ and $g : \o ^\b \hookrightarrow \bigcup _{i\in K}L_i \cap A$.
Let $i_0 =\max \{ i\in K : g[\o ^\b ] \cap L_i \cap A \neq \emptyset \}$. Then $F=g[\o ^\b ] \cap L_{i_0} \cap A$ is a final part of the linear order
$g[\o ^\b ] \cong \o ^\b $ and, since $\type (g[\o ^\b ]\setminus F)<\o ^\b $, by Fact \ref{T4400}(c) we have
$\type (F)=\o ^\b$ and, hence $\o ^\b \hookrightarrow L_{i_0}\cap A$ and $i_0 \geq j$. By (ii) of Claim \ref{T4437}
 we have $A\in (\Fin \times  \I _{\o ^\b }) ^+ $ and (\ref{EQ4405}) is proved. So (\ref{EQ4411}) is true.

By (\ref{EQ4405}) we have $\I _\L  = \Fin \times \I _{\o ^\b }$. Since $\la \o ^{\b +1 }, \in \ra \cong \L $, by Fact \ref{T4056}(a) we have $\I _{\o ^{\b +1}}\cong \I _\L$ and, hence,
\begin{equation}\label{EQ4412}
\I _{\o ^{\b +1 }}\cong \Fin  \times \I _{\o ^\b } .
\end{equation}
Let us assume that the statements (a) and (b) are true for $n$. By (\ref{EQ4412}) we have
$\I _{\o ^{\b +n+1 }}\cong \Fin  \times \I _{\o ^{\b +n} } \cong \Fin \times (\Fin ^n   \times \I _{\o ^\b }) = \Fin ^{n+1}   \times \I _{\o ^\b }$.
By Fact \ref{T4401} we have
$\la \P (\o ^{\b +n +1 }), \subset \ra
=\la (\I _{\o ^{\b +n+1 }})^+ , \subset \ra
\cong \la (\Fin ^{n+1} \times \I _{\o ^\b } )^+ , \subset  \ra$.

(c) follows from (a) and Fact \ref{T4042}(b) and (e).

(d) We use induction. For a proof of (d) for $n=1$
we show that the mapping
$
F: \la P(\o \times \o ^\b )/=_{\Fin \times \I _{\o ^\b }}, \leq _{\Fin \times \I _{\o ^\b }} \ra
\rightarrow \la \la P(\o ^\b )/=_{\I _{\o ^\b }} , \leq _{\I _{\o ^\b }} \ra ^\o /\equiv , \leq _\equiv\ra ,
$
given by $F([A]_{=_{\Fin \times \I _{\o ^\b }}} )=[\la [\pi _{\o ^\b}[A\cap L_i]]_{=_{\I _{\o ^\b }}} :i\in \o \ra]_\equiv$,
is an isomorphism.
\begin{cla} \rm \label{T4438}
For $A,B \subset \o \times \o ^\b$ we have: $A =_{\Fin \times \I _{\o ^\b }} B$ if and only if
\begin{equation}\label{EQ4408}
[\la [\pi _{\o ^\b}[A\cap L_i]]_{=_{\I _{\o ^\b }}} :i\in \o \ra]_\equiv =
[\la [\pi _{\o ^\b}[B\cap L_i]]_{=_{\I _{\o ^\b }}} :i\in \o \ra]_\equiv .
\end{equation}
\end{cla}
\dok
First, by (\ref{EQ4400}) we have
\begin{equation}\label{EQ4409}
A =_{\Fin \times \I _{\o ^\b }} B \Leftrightarrow
\exists j\in \o \;\; \forall i\geq j \;\; \pi _{\o ^\b } [(A\bigtriangleup B)\cap L_i] \in \I _{\o ^\b }   .
\end{equation}
On the other hand, (\ref{EQ4408}) holds iff there is $j\in \o$ such that for all $i\geq j$ we have
$\pi _{\o ^\b}[A\cap L_i]\bigtriangleup \pi _{\o ^\b}[B\cap L_i] \in \I _{\o ^\b }$, that is, since the restriction $\pi _{\o ^\b}\mid L_i$ is a bijection,
$(\pi _{\o ^\b}\mid L_i)[(A\cap L_i)\bigtriangleup (B\cap L_i)]= \pi _{\o ^\b}[(A\bigtriangleup B)\cap L_i)]\in \I _{\o ^\b }$.
\kdok
By Claim \ref{T4438}, $F$ is a well-defined injection.

For $[\la [X_i]_{=_{\I _{\o ^\b }}} :i\in \o \ra]_\equiv \in ( P(\o ^\b )/=_{\I _{\o ^\b }} ) ^\o /\equiv $
we have $F([A]_{=_{\Fin \times \I _{\o ^\b }}} )= [\la [X_i]_{=_{\I _{\o ^\b }}} :i\in \o \ra]_\equiv$, where $A=\bigcup _{i\in I}\{ i\}\times X_i $, so
$F$ is a surjection.

By (\ref{EQ4400}) we have $[A]_{=_{\Fin \times \I _{\o ^\b }}} \leq _{\Fin \times \I _{\o ^\b }} [B]_{=_{\Fin \times \I _{\o ^\b }}} $ iff
\begin{equation}\label{EQ4410}
\exists j\in \o \;\; \forall i\geq j \;\; \pi _{\o ^\b } [A\setminus B \cap L_i] \in \I _{\o ^\b }
\end{equation}
and $[\la [\pi _{\o ^\b}[A\cap L_i]]_{=_{\I _{\o ^\b }}} :i\in \o \ra]_\equiv \leq _\equiv
[\la [\pi _{\o ^\b}[B\cap L_i]]_{=_{\I _{\o ^\b }}} :i\in \o \ra]_\equiv$ iff there is $j\in \o$ such that for all $i\geq j$ we have
$\pi _{\o ^\b}[A\cap L_i]\setminus \pi _{\o ^\b}[B\cap L_i] \in \I _{\o ^\b }$, that is, since the restriction $\pi _{\o ^\b}\mid L_i$ is a bijection,
$(\pi _{\o ^\b}\mid L_i)[(A\cap L_i)\setminus (B\cap L_i)]= \pi _{\o ^\b}[A\setminus B\cap L_i)]\in \I _{\o ^\b }$.
Thus $F$ is an isomorphism.

Assuming that the statement is true for $n$, by (b) and (d) for $n=1$ we have
$P(\o ^{n+1} \times \o ^\b )/(\Fin ^{n+1} \times \I _{\o ^\b }) \cong  P(\o \times (\o ^n \times \o ^\b ))/(\Fin \times (\Fin ^n \times \I _{\o ^\b }))
                                                                \cong  P(\o \times \o ^{\b +n})/(\Fin \times \I _{\o ^{\b +n} }))
                                                                \cong  \rp (P(\o ^{\b +n} )/ \I _{\o ^{\b +n} })
\\                                                                \cong  \rp (P(\o ^n \times \o^\b  )/ (\Fin ^n \times \I _{\o ^\b })
                                                                \cong  \rp (\rp ^n (P(\o ^\b )/ \I _{\o ^\b }))
                                                                \cong  \rp ^{n+1} (P(\o ^\b )/ \I _{\o ^\b }) $.

(e) follows from (c) and (d).
\kdok
\noindent
For $n\in \N$, let the ideal $\Fin ^n $ on the set $\o ^n = \o \times ( \o \times \dots \times (\o \times \o ) \dots )$ ($n$-many factors) be defined by:
$\Fin ^n = \Fin \times ( \Fin \times \dots \times (\Fin \times \Fin ) \dots )$ ($n$-many factors). Then, by Lemma \ref{T4403} we have
\begin{cor}\rm\label{T4404}
For each $n\in \N$ we have:

(a) $\la \P (\o ^n ), \subset \ra \cong \la P(\o ^n)\setminus \Fin ^n  , \subset  \ra$ and  $\I _{\o ^n }\cong \Fin ^n $;

(b) $\sq (\P (\o ^n ), \subset ) \cong (\rp ^{n-1} (P(\o )/\Fin  ))^+ $.
\end{cor}
\begin{lem}\rm\label{T4409}
$\la \P (\g +k ),\subset \ra \cong \la  \P (\g ),\subset \ra$, for each limit ordinal $\g$ and each $k\in \N$.
\end{lem}
\dok
%Clearly we have $\g +k = \g \cup \{ \g , \g +1, \dots , \g +k-1\}$ and
First we prove
$\P (\g +k )= \{ C \cup \{ \g , \g +1, \dots , \g +k-1\} : C\in \P (\g ) \} $.
The inclusion ``$\supset$" is evident. If $A\in \P (\g +k )$ and $f:\g +k \hookrightarrow \g +k$, where $A=f[\g +k]$, then,
since $f$ is an increasing function, we have $f(\b )\geq \b$, for each $\b \in \g +k$, which implies $f(\g +i) =\g +i$, for $i<k$,
and, hence, $C= f[\g ] \in \P (\g )$ and $A=C \cup \{ \g , \g +1, \dots , \g +k-1\}$.

Now it is easy to show that the mapping $F: \la \P (\g ),\subset \ra \rightarrow \la \P (\g +k ),\subset\ra$,
given by $F(C)=C \cup \{ \g , \g +1, \dots , \g +k-1\} $, is an isomorphism.
\hfill $\Box$
\begin{lem} \rm \label{T4436}
Let $\delta , \delta ' >0$  be countable ordinals. Then

(a) The ordinal $\o ^\delta$ is an $\o$-sum of elements of $\H$ satisfying (\ref{EQ4100});

(b) $\delta \geq \delta ' \Rightarrow \o ^\delta +\o ^{\delta '} \not\in \H$.
\end{lem}
\dok
(a) By Fact \ref{T4400} we have $\o ^\delta \in \H$ and $\o ^\delta $ can not be an $\o ^*$-sum (since it is a well ordering) so it is an
$\o$-sum of elements of $\H$ satisfying (\ref{EQ4100}).

(b) Suppose that $\o ^\delta +\o ^{\delta '} \in \H$. Then, by Fact \ref{T4400}, $\o ^\delta +\o ^{\delta '}=\o ^{\delta ''}$, for some
ordinal $\delta ''$ and, clearly $\o ^\delta \leq \o ^{\delta ''}$. Now, $\o ^\delta =\o ^{\delta ''} $ is impossible, since $\o ^\delta$ can not be isomorphic
to its proper initial segment and, hence, $\o ^\delta <\o ^{\delta ''} $, which implies that $\o ^{\delta '} <\o ^{\delta ''} $ as well.
But this is impossible by Fact \ref{T4400}(b).
\kdok
\noindent
{\bf Proof of Theorem \ref{T4407}.}
By Lemma \ref{T4409} we can assume that $k=0$.
So, we have $\a = \o ^{\g _n +r_n } + \dots + \o ^{\g _n +r_n } + \dots + \o ^{ \g _0 +r_0 }+ \dots +\o ^{ \g _0 +r_0 }=\sum _{j< \sum _{i=0}^n s_i}L_j$.
By Lemma \ref{T4436}(a) for each $j$ the order $L_j\in \H$ and it is an $\o$-sum of elements of $\H$ satisfying (\ref{EQ4100}).
By Lemma \ref{T4436}(b), $L_j + L_{j+1}\not\in \H $ so, by Fact \ref{T4408}(b),
%\begin{equation}\label{EQ4412}\textstyle
$
\la \P (\a ), \subset \ra
\cong \prod _{j< \sum _{i=0}^n s_i} \la \P (L_j) , \subset \ra
= \prod _{i=0}^n \la \P (\o ^{\g _i + r_i}) , \subset \ra ^{s_i}
$,
%\end{equation}
which, with Fact \ref{T4042}(b) and (c) and Lemma \ref{T4403}(e), gives
$
\sq \la \P (\a ), \subset \ra   \cong  \textstyle \sq \prod _{i=0}^n \la \P (\o ^{\g _i + r_i}) , \subset \ra ^{s_i}$
$                                \cong  \textstyle \prod _{i=0}^n (\sq \la \P (\o ^{\g _i + r_i}) , \subset \ra )^{s_i}
                               \cong  \textstyle \prod _{i=0}^n ( ( \rp ^{r_i}( P(\o ^{\g _i} )/ \I _{\o ^{\g _i} }))^+ )^{s_i}. %\hfill \Box
$
\hfill $\Box$
\begin{cor} \rm \label{T4410}
$\sq \la \P (\o n),\subset \ra \cong ((P(\o )/\Fin )^+ )^n$, for each $n\in \N$.
%\begin{equation}\label{EQ4415}\textstyle
%$\sq \la \P (\; \underbrace{\o + \dots + \o}_n\; ),\subset \ra \cong \underbrace{(P(\o )/\Fin )^+ \times \dots \times (P(\o )/\Fin )^+}_n$.
%\end{equation}
\end{cor}
\section{Forcing with the quotient $P(\o ^\g )/\I _{\o ^\g}$}\label{S4}
By Theorem \ref{T4407}, the poset $\la \P (\alpha ), \subset \ra$ is forcing equivalent to a forcing product
of iterated reduced products of Boolean algebras of the form $P(\o ^\g )/\I _{\o ^\g}$.
In this section we consider such algebras and assume that $\g \geq \o$ is a countable limit ordinal, $\la \d _n : n\in \o  \ra$ a fixed increasing cofinal sequence in $\g \setminus \{ 0\}$
and $\L =\la L, < \ra  $ $=\sum _{n\in \o}\la L_n , <_n\ra$, where $\la L_n , <_n \ra\cong \la \o ^{\d _n} , \in \ra$, for $n\in \o$, and
$L_m \cap L_n=\emptyset $, for $m\neq n$. For $A\subset L$ and $m\in \o$ let
$S^m_A=\{ n\in \o : \type (A\cap L_n )\geq \o ^{\d _m }\}$ and $\supp A=\{n\in \o : A\cap L_n \neq \emptyset \}$.

The ideal $\I _\L=\{ A\subset L: \L \not\hookrightarrow A \}$ will be denoted by $\I$ and,
if $G\subset P(\o )$ is an ultrafilter, $\I _G=\{ A\subset L : \exists I \in \I \; \supp (A\setminus I)\not\in G \}$.
$\Gamma$ (resp.\ $\Gamma _1$) will be the canonical name for a $\la [\o ]^\o , \subset ^* \ra$-generic (resp.\ $(P(\o )/\Fin )^+$-generic) filter over the ground model $V$
and $q: P(\o ) \rightarrow P(\o )/\Fin $ the quotient mapping.

The aim of this section is to prove the following statement. It follows from Propositions \ref{T4425} and \ref{T4428} given at the end of the section.
\begin{te}  \rm  \label{T4420}
For each countable limit ordinal $\g$ we have:

(a) The partial orders $\langle \P (\o ^\g) , \subset \rangle$ and $(P(\o ^\g )/\I _{\o ^\g})^+$
are forcing equivalent to the two-step iteration
$
(P(\o )/\Fin )^+ \ast (\check{P(L)}/\check{\I } _{\check{q}^{-1}[\Gamma _1]})^+ ;
$

(b) $[\o ]\Vdash $ ``$(\check{P(L)}/\check{\I } _{\check{q}^{-1}[\Gamma _1]})^+$ is an $\o _1$-closed, separative and atomless poset".
\end{te}
\begin{fac}\rm\label{T4413}
Let $f:\o \rightarrow \o$ be an increasing function. Then

(a) $\o ^{\d _{f(0)}} + \o ^{\d _{f(1)}}+\dots +\o ^{\d _{f(m)}} =\o ^{\d _{f(m)}}$, for each $m\in \o$;

(b) $\sum _{n\in \o}\o ^{\d _{f(n)}}= \o ^\g$;

(c) $\L \cong \o ^\g$.
\end{fac}
\dok
We prove (a) by induction. Assuming that (a) is true for $m\in \o$ we have
$\o ^{\d _{f(0)}} + \dots +\o ^{\d _{f(m+1)}}
%=(\o ^{\d _{f(0)}} +\dots +\o ^{\d _{f(m)}}) +  \o ^{\d _{f(m+1)}}
=\o ^{\d _{f(m)}}+ \o ^{\d _{f(m+1)}}
=\o ^{\d _{f(m)}}\cdot 1 + \o ^{\d _{f(m)} + (\d _{f(m+1)}- \d _{f(m)})}
=\o ^{\d _{f(m)}}(1 +\o ^{\d _{f(m+1)}- \d _{f(m)}})
=\o ^{\d _{f(m)}}\o ^{\d _{f(m+1)}- \d _{f(m)}}
=\o ^{\d _{f(m+1)}} $.

(b) By (a) and basic properties of ordinal arithmetic we have
$\sum _{n\in \o}\o ^{\d _{f(n)}}= \sup \{ \sum _{n\leq m}\o ^{\d _{f(n)}} : m\in \o  \}=\sup \{ \o ^{\d _{f(m)}} : m\in \o \}=\o ^\g$.

(c) By (b) we have $\L \cong \sum _{n\in \o}\o ^{\d _n}= \o ^\g$.
\hfill $\Box$
\begin{lem}\rm\label{T4422}
For $A\subset L$ and $m\in \o$ we have:

(a) $S^m_A \subset \supp A \setminus m$;

(b) $ m_1 <m_2 \Rightarrow S^{m_1}_A \supset S^{m_2}_A $;

(c) $ A\subset B \Rightarrow S^{m}_A \subset S^{m}_B $;

(d) $A\in P(L)\setminus \I $ iff $S^m _A \in [\o ]^\o$, for each $m\in \o$;

(e) $A\in \I $ iff $S^m _A=\emptyset$, for some $m\in \o$;

(f) $|\supp (A)|<\o \Rightarrow A\in \I$;

(g) $S^m_{\bigcup _{k<l}A_k}=\bigcup _{k<l}S^m_{A_k}$;

(h) $A\subset _{\I }B$ iff $S^m _{A\setminus B}=\emptyset$, for some $m\in \o$.
\end{lem}
\dok
(a), (b), (c) and (f) are evident and  (h) follows from (e).

(d) By Fact \ref{T4408}, $A\in P(L)\setminus \I $ iff for each $m\in \o$ we have: for each $n\in \o$ there is finite
$K\subset \o \setminus n$ such that $L_m \cong \o ^{\d _m}\hookrightarrow \bigcup _{i\in K}A\cap L_j$, but, by Fact \ref{T4400}, $\o ^{\d _m}$ is
an indivisible structure and, hence, this holds iff there is $k\geq n$ such that $\o ^{\d _m}\hookrightarrow A\cap L_k$, that is $k\in S^m_A$.

(e) By (c), if $S^m_A=\emptyset$ for some $m\in \o$, then $A\in \I$. On the other hand, if $A\in \I$, then, by (c) again,
there are $k,l\in \o$ such that $S^k_A \subset l$ and, by (a) and (b), for $m\geq l,k$ we have
$S^m_A \subset S^k_A \setminus m\subset S^k_A \setminus l=\emptyset$.

(g) If $n\in S^m_{\bigcup _{k<l}A_k}$, then $\o ^{\d _m}\hookrightarrow \bigcup _{k<l}A_k \cap L_n$ and, since
$\I _{\o ^{\d _m}} =\{ I\subset \o ^{\d _m}: \o ^{\d _m}\not\hookrightarrow I \}$ is an ideal, there is $k<l$ such that
$\o ^{\d _m}\hookrightarrow A_k \cap L_n$, that is $k\in S^m_{A_k}$. On the other hand, by (c) we have $S^m_{A_k}\subset S^m_{\bigcup _{k<l}A_k}$, for each $k<l$.
\hfill $\Box$
\begin{lem}\rm\label{T4427}
If $G\subset P(\o )$ is an ultrafilter, then

(a) $\I _G=\{ A\subset L : \exists I \in \I \; \supp (A\setminus I)\not\in G \}$ is an ideal and $\I\subset \I _G$;

(b) $\sm \la P(L)\setminus \I _G , \subset _\I \ra =\la P(L)\setminus \I _G , \subset _{\I _G} \ra $.
%where instead $\subset _\I \cap (P(L)\setminus \I _G)^2$ we write $\subset _\I$.
\end{lem}
\dok
(a) If $A_1,A_2\in \I _G$ and $\supp (A_1\setminus I_1) , \supp (A_2\setminus I_2)\not \in G$, where $I_1 , I_2 \in \I$, then, since
$(A_1 \cup A_2)\setminus (I_1 \cup I_2)\subset (A_1\setminus I_1) \cup (A_2\setminus I_2)$ and $\supp (X\cup Y)=\supp (X)\cup \supp (Y)$, we have
$\supp ((A_1 \cup A_2)\setminus (I_1 \cup I_2)) \subset \supp (A_1\setminus I_1) \cup \supp (A_2\setminus I_2)\not \in G$ and, since $I_1 \cup I_2\in \I$,
we have $A_1 \cup A_2\in \I _G$.

(b) Let $A\subset _{\I _G}B$ and $C\in P(L)\setminus \I _G$, where $C\subset _\I A$. Then, since $C= (C\setminus A)\cup (C\cap A\setminus B)\cup (C\cap A \cap B)$,
$A\setminus B \in \I _G$ and $C\setminus A\in \I \subset \I _G$,
we have $D=C\cap A \cap B\in P(L)\setminus \I _G$ and $D \subset _\I C,B$. Thus $A\subset _\I ^* B$.
Conversely, suppose that $A\subset _\I ^* B$ and $A\setminus B \not\in \I _G$. Then, for $C=A\setminus B$ there is $D\in P(L)\setminus \I _G$
such that $D\subset _\I A\setminus B$ and $D\subset _\I B$, which implies $D\in \I$. A contradiction.
\kdok

\noindent
We remind the reader that, if $\la \P \leq _\P , 1_\P \ra$ and  $\la \Q \leq _\Q , 1_\Q \ra$ are pre-orders, then a mapping
$f: \P \rightarrow \Q$ is a {\it complete embedding}, in notation $f:\P \hookrightarrow _c \Q$  iff

(ce1) $p_1 \leq _\P p_2 \Rightarrow f(p _1)\leq _\Q f(p_2)$,

(ce2) $p_1 \perp _\P p_2 \Leftrightarrow f(p _1)\perp _\Q f(p_2)$,

(ce3) $\forall q\in \Q \; \exists p\in \P \; \forall p' \leq _\P p \;\; f(p') \not\perp _\Q q$.

\noindent
Then, for $q\in \Q $ the set $\red (q) =\{ p\in \P : \forall p' \leq _\P p \; f(p') \not\perp _\Q q \}$ is the set of
{\it reductions} of $q$ to $\P$. The following fact is folklore (see \cite{Kun}).
\begin{fac}\rm\label{T4432}
If $f:\P \hookrightarrow _c \Q$, then $\Q$ is forcing equivalent to the two step iteration $\P \ast \la \pi , \leq _\pi , \check{1_\Q}\ra $, where
$1_\P \Vdash _\P \pi \subset \check{\Q}$ and for each $p\in \P$ and $q,q_1,q_2\in \Q$

(a) $p\Vdash \check{q}\in \pi$ iff $p\in \red (q)$;

(b) $p\Vdash \check{q_1}\leq _\pi \check{q_2}$ iff $q_1 \leq _\Q q_2$ and $p\in \red (q_1)$.
\end{fac}

\begin{prop}\rm\label{T4425}
The following pre-orders are  forcing equivalent:

1. $\la \P (\o ^\g), \subset \ra $,

2. $(P(\o ^\g )/\I _{\o ^\g })^+$,

3. $\la P(L)\setminus \I , \subset _\I \ra$,

4. $\la [\o ]^\o ,\subset ^*\ra \ast \la \check{P(L)}\setminus \check{\I } _{\Gamma } , \subset _{\check{\I }_\Gamma }  \ra$,

5. $(P(\o )/\Fin )^+ \ast (\check{P(L)}/\check{\I } _{\check{q}^{-1}[\Gamma _1]})^+$.
\end{prop}
\dok
By Facts \ref{T4401}, \ref{T4042}(a) and (e) the posets 1 and 2 are forcing equivalent.
By Facts \ref{T4401}, \ref{T4056}(b) and \ref{T4042}(a),(d) the poset
$\la \P (\o ^\g), \subset \ra
= \la P (\o ^\g)\setminus \I _{\o ^\g }, \subset \ra$
is isomorphic to the poset $\la P (L)\setminus \I , \subset \ra$, forcing equivalent to
$\la P (L)\setminus \I , \subset _{\I }\ra$.
The forcing equivalence of the posets  4 and 5 is evident - note that $G_1$ is a $(P(\o )/\Fin )^+$-generic filter iff $G=q^{-1}[G_1]$ is an
$\la [\o ]^\o , \subset ^* \ra$-generic filter over $V$ and that
$\sq \la P(L)\setminus \I _G , \subset _{\I _G}\ra =(P(L)/\I _G )^+$.

Thus the forcing equivalence of the posets  3 and 4 remains to be proved.
\begin{cla}\rm\label{T4433}
The mapping $f: \la [\o ]^\o ,\subset ^*\ra \rightarrow \la P(L)\setminus \I , \subset _\I \ra$ defined by
$f(S)=\bigcup _{n\in S}L_n$ is a complete embedding. In addition, $\t (\la P(L)\setminus \I , \subset _\I \ra)\leq \t$.
\end{cla}
\dok
By Fact \ref{T4413}(b) and (c), for $S\in [\o ]^\o$ we have
$f(S)\cong \sum _{n\in S}\o ^{\d _n}=\o ^\g \cong \L$, thus  $f(S)\in P(L)\setminus \I$. Let $S,T\in [\o ]^\o$.

(ce1) If $S\subset ^* T$, then $|\supp (f(S)\setminus f(T))|=|\supp (\bigcup _{n\in S\setminus T }L_n)|=|S\setminus T|<\o$ and,
by Lemma \ref{T4422}(f), $f(S)\setminus f(T)\in \I$, that is $f(S)\subset _\I f(T)$.

(ce2) If $S\perp T$, then $|\supp (f(S)\cap f(T))|=|\supp (\bigcup _{n\in S\cap T }L_n)|=|S\cap T|<\o$ and,
by Lemma \ref{T4422}(f), $f(S)\cap f(T)\in \I$, that is $f(S)\perp _\I f(T)$.
If $S\not\perp T$, then $S\cap T\in [\o ]^\o $ and $f(S)\cap f(T)=\bigcup _{n\in S\cap T }L_n= f(S\cap T)\in P(L)\setminus \I $ and,
hence,  $f(S)\not\perp _\I f(T)$.

(ce3) First we show that for $S\in [\o ]^\o$ and $A\in P(L)\setminus \I $ we have
\begin{equation}\label{EQ4438}
S\in \red (A) \mbox{ iff }S\subset ^* S^m _A, \mbox{ for each }m\in \o.
\end{equation}
Suppose that $S\in \red (A)$ and that $T=S\setminus S^m_A\in [\o]^\o$, for some $m\in \o$. Then there is
$B\in P(L)\setminus \I$ such that $B\subset _\I f(T), A$. Now we use Lemma \ref{T4422}. By (h), there are $m_1 , m_2 \in \o$
such that $S^{m_1}_{B\setminus f(T)}=S^{m_2}_{B\setminus A}=\emptyset$.  By (b), for $m^*=\max \{m,m_1,m_2\}$ we have
$S^{m^*}_{B\setminus f(T)}=S^{m^*}_{B\setminus A}=\emptyset$ and, by (g),
$S^{m^*}_B=S^{m^*}_{(B\cap A \cap f(T))\cup (B\setminus f(T))\cup (B\setminus A)}=
S^{m^*}_{B\cap A \cap f(T)}\cup S^{m^*}_{B\setminus f(T)}\cup S^{m^*}_{B\setminus A}=S^{m^*}_{B\cap A \cap f(T)}$.
But, by (a), (b) and (c), $S^{m^*}_{B\cap A \cap f(T)}\subset S^{m^*}_{f(T)}\cap S^{m^*}_{A}\subset T\cap S^m_A=\emptyset$, that is $S^{m^*}_B=\emptyset$,
which, by (e), implies $B\in \I$. A contradiction.

Let $S\subset ^* S^m _A$,  for each $m\in \o$, and let $[\o ]^\o \ni T \subset ^* S$.
In order to find a set $B\in P(L)\setminus \I$ such that $B\subset _\I f(T), A$
by recursion we construct a sequence $\la n_k : k\in \o \ra$ such that for each $k\in \o$ we have:
(i) $n_k \in T$, (ii) $n_k < n_{k+1}$ and (iii) $\type (A\cap L_{n_k})\geq \o ^{\d _k}$.
If a sequence $\la n_0, \dots , n_k \ra$ satisfies (i)-(iii), then for $n_{k+1}=\min (T\cap S^{k+1}_A\setminus (n_k +1))$ we have
$n_k <n_{k+1}$ and $\type (A\cap L_{n_{k+1}})\geq \o ^{\d _{k+1}}$ thus the recursion works.
Now $B=\bigcup _{k\in \o}A\cap L_{n_k}\subset A$, by (i) we have $B\subset f(T)$
and, by (iii), for each $k\in \o $ we have $S^k_B \supset \{ n_k , n_{k+1}, \dots \}$. Thus, by Lemma \ref{T4422}(d) we have
$B\in P(L)\setminus \I$ and (\ref{EQ4438}) is proved.

Now we check (ce3). If $A\in P(L)\setminus \I $ then, by Lemma \ref{T4422} $\{ S^m_A : m\in \o \}$ is a subfamily of $[\o ]^\o$
having the strong finite intersection property and, hence, it has a pseudointersection $S\in [\o]^\o$. By  (\ref{EQ4438}), $S$ is a reduction of $A$ to $[\o ]^\o$.

If $\la T_\a : \a <\t \ra$ is a tower in $\la [\o ]^\o ,\subset ^*\ra$, then, by (ce1), $\la f(T_\a) : \a <\t \ra$
is a $\subset _\I$-decreasing sequence in $P(L)\setminus \I$. Suppose that $A\in P(L)\setminus \I$ and that for each $\a <\t$ we have
$A\subset _\I f(T_\a)$, which, by Lemma \ref{T4422}(h), gives $m_\a \in \o$ such that $S^{m_\a}_{A\setminus f(T_\a)}=\emptyset$ and,
by Lemma \ref{T4422}(g), $S^{m_\a}_A=S^{m_\a}_{A\cap f(T_\a)}\subset T_\a$.
Let $S\in \red (A)$. Then, by (\ref{EQ4438}), $S\subset ^* S^{m_\a}_A\subset T_\a$, for each $\a <\t$. A contradiction.
\kdok
By the previous claim, Fact \ref{T4432} and (\ref{EQ4438}), the pre-order $\la P(L)\setminus \I , \subset _\I \ra$ is forcing equivalent to the iteration
$\la [\o ]^\o ,\subset ^*\ra  \ast \la \pi , \leq _\pi , \check{L} \ra $,
where $\o \Vdash \pi \subset (P(L)\setminus \I )\check{\;} $
%$\pi =\{ \la \check{A}, S \ra : A\in P (L)\setminus \I  \land S\in [\o ]^\o \land \forall m\in \o \; S\subset ^* S^m _A \}$,
%$\leq _\pi = \{ \la\la A,B\ra \check{\;} , S\ra : A,B \in P (L)\setminus \I \land S\in [\o ]^\o \land
%A\subset _\I B \land \forall m\in \o \; S\subset ^* S^m _A\}$
%$1_\pi =L\check{\;}$.
 and for each $S\in [\o ]^\o $ and $A,B \in P(L)\setminus \I $ we have
\begin{eqnarray}
S\Vdash \check{A}\in \pi             & \Leftrightarrow & \forall m\in \o \;\; S\subset ^* S^m _A ;\label{EQ4439}\\
S\Vdash \check{A}\leq _\pi \check{B} & \Leftrightarrow & S\Vdash \check{A}\in \pi \;\land \; A\subset _\I B \label{EQ4440}.
\end{eqnarray}
\begin{cla}\rm\label{T4434}
(a) $\o \Vdash ``\check{\I }$ and $\check{\I}_\Gamma$ are ideals in $P(\check{L})=\check{P(L)}"$;

(b) $\o \Vdash \pi = \check{P(L)}\setminus \check{\I}_\Gamma$;

(c) $\o \Vdash \leq _\pi = \subset _{\check{\I}}\cap (\pi \times \pi)$.
\end{cla}
\dok
Let $G$ be an $\la [\o ]^\o ,\subset ^*\ra$-generic filter over $V$.

(a) Since the forcing $\la [\o ]^\o , \subset ^* \ra$ is $\o$-distributive, in $V[G]$ we have  $P^{V[G]}(L)=P^V(L)$ and, for the same reason,
$\I$ remains to be an ideal in $P^{V[G]}(L)$. By Lemma \ref{T4427}(a), the set $\I _G=\{ A\subset L : \exists I \in \I \; \supp (A\setminus I)\not\in G \}$ is an ideal in $P^{V[G]}(L)$.

(b) We show that
$\pi _G =\{ A\in P(L)\setminus \I : \forall I\in \I \; \supp (A\setminus I)\in G \}$.
Let $A\in P(L)\setminus \I$.
If $A\in \pi _G$, then there is $S\in G $ such that $S\Vdash \check{A}\in \pi$. For $I\in \I$ we have
$A\cap I\in \I$ and, by Lemma \ref{T4422}(e), there is $m^*\in \o$ such that $S^{m^*}_{A\cap I}=\emptyset$.
Thus, by (\ref{EQ4439}) and Lemma \ref{T4422}(g) we have
$S\subset ^* S^{m^*}_A = S^{m^*}_{A\cap I}\cup S^{m^*}_{A\setminus I}=S^{m^*}_{A\setminus I}\subset \supp(A\setminus I)$, which implies that
$\supp(A\setminus I)\in G$. So $A\in P(L)\setminus \I _G$.

If $A\not\in \pi _G$, then there is $S\in G $ such that $S\Vdash \neg \check{A}\in \pi$. Suppose that
$|S \cap S^m_A|=\o$, for each $m\in \o$. Then, by Lemma \ref{T4422}(b), $S \cap S^m_A$, $m\in \o$, would be a decreasing sequence in $[\o ]^\o$ and,
hence, there would be $T\in [\o ]^\o$ such that $T\subset ^* S \cap S^m_A$, for each $m\in \o$, which, by (\ref{EQ4439}), implies
$T\Vdash \check{A}\in \pi$. But this is impossible since $T\subset ^* S$ and $S\Vdash \neg \check{A}\in \pi$.
Thus $|S \cap S^{m^*}_A|<\o$, for some $m^* \in \o$.
Let $I=\bigcup _{n\in S}A\cap L_n$. By Lemma \ref{T4422}(c) we have $S^{m^*}_I \subset S \cap S^{m^*}_A$ and, hence $|S^{m^*}_I|<\o$, which,
by Lemma \ref{T4422}(d), implies $I\in \I$.
Since $\supp (A\setminus I)\cap S=\emptyset$ and $S\in G$, we have $\supp (A\setminus I)\not\in G$. So $A\not\in P(L)\setminus \I _G$.

(c) We show that
$(\leq _\pi )_G= \subset _\I \cap (P(L)\setminus \I _G)^2$.
Let $A,B \in \pi _G$ and let $S\in G $ where $S\Vdash \check{A},\check{B}\in \pi$.
If $A (\leq _\pi )_G B$, then there is $T\in G$ such that $T\Vdash \check{A} \leq _\pi \check{B}$, which, by (\ref{EQ4440}),
implies $A\subset _\I B$. If $A\subset _\I B$, then, since $S\Vdash \check{A}\in \pi$, by (\ref{EQ4440}) we have
$S\Vdash \check{A} \leq _\pi \check{B}$ and, hence,  $A (\leq _\pi )_G B$.
\kdok
Thus, the pre-order $\la P(L)\setminus \I , \subset _\I \ra$ is forcing equivalent to the two-step iteration
$\la [\o ]^\o ,\subset ^*\ra \ast
\la \check{P(L)}\setminus \check{\I } _{\Gamma } , \check{\subset _\I } \ra$ and, by
Lemma \ref{T4427}(b) applied in $V[G]$, to the iteration $\la [\o ]^\o ,\subset ^*\ra \ast
\la \check{P(L)}\setminus \check{\I } _{\Gamma } , \subset _{\I _\Gamma }  \ra$.
\hfill $\Box$
\begin{prop}\rm\label{T4428}
According to the notation of Proposition \ref{T4425} we have

(a) $\o \Vdash ``\la \check{P(L)}\setminus \check{\I } _{\Gamma } , \subset _{\check{\I } _\Gamma }  \ra $ is a separative, $\o _1$-closed and atomless pre-order".

(b) $[\o ]\Vdash ``(\check{P(L)}/\check{\I } _{\check{q}^{-1}[\Gamma _1 ]})^+$ is a separative, $\o _1$-closed and atomless poset."
\end{prop}
\dok
(a)
The separativity follows from Fact \ref{T4042}(d) and we prove $\o _1$-closure.
We easily show that for $S\in [\o ]^\o$ and $A,B\in P(L)$ satisfying $S\Vdash \check{A}, \check{B}\in \pi$
we have:
\begin{equation}\label{EQ4441}\textstyle
S\Vdash \check{A}\subset _{\check{\I }_\Gamma }\check{B} \;\Leftrightarrow\;
\forall T\subset ^* S \;\; \exists I\in \I \;\; |T\setminus \supp (A\setminus B)\setminus I))|=\o.
\end{equation}
Since the forcing $\la [\o]^\o  , \subset ^* \ra$ is $\o$-distributive we have
$\o \Vdash \check{P(L)}^{\check{\o }}=((P(L)^\o )^V )\check{\;}$ and, clearly,
$\o \Vdash \pi \subset \check{P(L)}$.
So, assuming that $\la A_n :n\in \o\ra \in P(L)^\o$, $S\in [\o ]^\o$ and
$S \Vdash \forall n\in \check{\o } \; (\check{A_n}\in \pi \land \forall m\geq n \;\check{A_m}\subset _{\check{\I }_\Gamma} \check{A_n})$
that is, by (\ref{EQ4439}) and (\ref{EQ4441}),
\begin{equation}\label{EQ4442}
\forall m,n \in \o \;\; S\subset ^* S^m_{A_n} ,
\end{equation}
\begin{equation}\label{EQ4443}
\forall R\subset ^* S \;\; \exists I\in \I \;\; |R\setminus \supp ((A_{n+1}\setminus A_n) \setminus I)|=\o  ,
\end{equation}
it is sufficient to find $T\in [\o ]^\o$ and $A\in P(L)$ such that $T\subset ^* S$,
$T\Vdash \check{A}\in \pi $ and $T\Vdash \check{A}\subset _{\check{\I }_\Gamma} \check{A_n}$, for all $n\in \o $.
\begin{cla}\rm\label{T4435}
For $r\in \o$, let $\textstyle S_r =S\cap\bigcap _{m,n\leq r}S^m _{A_n} \mbox{ and } B_r =A_r \cap \bigcup _{k\in S_r}L_k$. Then

(a) $B_r\in P(L)\setminus \I$;

(b) $B_{r+1}\subset _\I B_r$.
\end{cla}
\dok
(a) If $m\in \o$, then $k\in S^m _{B_r}$ iff $k\in S_r$ and $\o ^{\d _m}\hookrightarrow B_r \cap L_k=A_r \cap L_k$; thus $S^m_{B_r}=S_r\cap S^m_{A_r}$ and,
by (\ref{EQ4442}), $|S^m_{B_r}|=\o$. Now, by Lemma \ref{T4422}(d), $B_r\in P(L)\setminus \I$.

(b) Suppose that $B_{r+1}\not\subset _\I B_r $. Then, since $S_{r+1}\subset S_r$, we would have
$C= B_{r+1}\setminus B_r =(A_{r+1}\setminus A_r) \cap \bigcup _{k\in S_{r+1}}L_k \in P(L)\setminus \I$ and, by Lemma \ref{T4422}(b) and (d),
there would be $R\in [\o ]^\o$ such that
$R\subset ^* S^m_C $,
for all $m \in \o$.
Since $S^m_C \subset \supp (C)\subset S_{r+1}\subset S$ we would have $R\subset ^* S$ and, by (\ref{EQ4443}), there would be $I\in \I$ such that
$R\not\subset ^* \supp ((A_{r+1}\setminus A_r) \setminus I)) $.
Since $(A_{r+1}\setminus A_r) \cap I\in \I$, by Lemma \ref{T4422}(e) there is $m^*\in \o$ such that
$S^{m^*} _{(A_{r+1}\setminus A_r) \cap I}$ $=\emptyset$ and, by Lemma \ref{T4422}(g),
$S^{m^*} _{A_{r+1}\setminus A_r}
%=S^{m^*} _{(A_{r+1}\setminus A_r) \cap I} \cup S^{m^*} _{(A_{r+1}\setminus A_r) \setminus I}
=S^{m^*} _{(A_{r+1}\setminus A_r) \setminus I}\subset \supp ((A_{r+1}\setminus A_r) \setminus I)$.
But, by Lemma \ref{T4422}(c),
$R\subset ^* S^{m^*} _C \subset S^{m^*} _{A_{r+1}\setminus A_r}$  thus $R\subset ^*\supp ((A_{r+1}\setminus A_r) \setminus I)$. A contradiction.
\kdok
By Theorem \ref{T4121} the pre-order $\la P(L)\setminus \I , \subset _\I \ra$ is $\o _1$-closed so, by Claim \ref{T4435}, there is $A\in P(L)\setminus \I $ such that
\begin{equation}\label{EQ4447}
\forall n\in \o\;\;A\subset _\I B_n \subset A_n  .
\end{equation}
By Lemma \ref{T4422}(b) and (d) there is $T\in [\o ]^\o$ such that
\begin{equation}\label{EQ4448}
\forall m \in \o \;\; T\subset ^* S^m_A .
\end{equation}
By (\ref{EQ4447}) we have $A\setminus B_n \in \I$, by Lemma \ref{T4422}(e) there is $m^*\in \o$ such that
$S^{m^*} _{A\setminus B_n} =\emptyset$ and, by Lemma \ref{T4422}(g) we have
$S^{m^*} _A
=S^{m^*} _{A\cap B_n} \cup S^{m^*} _{A \setminus B_n}
=S^{m^*} _{A\cap B_n}
\subset S^{m^*} _{B_n}
\subset \supp (B_n)
\subset S_n
\subset S$. By  (\ref{EQ4448}) we have $T\subset ^* S^{m^*} _A$ and, hence, $T\subset ^*S$.
By (\ref{EQ4448}) and (\ref{EQ4439}) we have $T\Vdash \check{A}\in \pi$. By (\ref{EQ4447}), for each $n\in \o$ we have
$A \subset _\I A_n$ and, hence, $T\Vdash \check{A}\subset _{\check{\I }_\Gamma} \check{A_n}$.

Taking an $\la [\o ]^\o , \subset ^*\ra$-generic filter $G$ we prove that the pre-order
$\la \pi _G , \subset _{\check{\I } _G }  \ra $ is atomless. If $A\in \pi _G$, then, by (\ref{EQ4439}), there is $S\in G$
such that $S \subset ^* S^m_A$, for each $m\in \o$. By Lemma \ref{T4422}(b) we have $S^0_A \supset S^1_A \supset \dots$ and, clearly,
$\bigcap _{m\in \o} S^m_A =\emptyset$. W.l.o.g.\ suppose that $S\subset S^0_A$. Then
$S=\bigcup _{m\in \o}S \cap (S^m_A \setminus S^{m+1}_A)$ and, for $n\in S \cap (S^m_A \setminus S^{m+1}_A)$ there is
$\f _n : \o ^{\d _m}\hookrightarrow A\cap L_n$. Let $\f _n [\o ^{\d _m}]=B_n \dot{\cup } C_n$, where
$B_n,C_n \cong \o ^{\d _m}$ and let $B=\bigcup _{m\in S}B_n$ and  $C=\bigcup _{m\in S}C_n$. Then
$S^m _B =S\cap S^m _A$ and, hence, $S \subset ^* S^m _B$, for all $m\in \o$, which implies $B\in \pi _G$ and, similarly,
$C\in \pi _G$. Since $B,C\subset A$ we have $B,C \subset _{\I _G}A$ and $B\cap C=\emptyset$ implies that $B$ and $C$ are
$\subset _{\check{\I } _G}$-incompatible.

The proof of (b) is similar to the proof of (a). Note that $\la P(L)\setminus \I _{G}, \subset _{\I _G}\ra$
is $\o _1$-closed (atomless) iff $(P(L)/ \I _{G})^+$ is $\o _1$-closed (atomless).
\hfill $\Box$
\section{Forcing with $\la \P(\a ), \subset \ra $}\label{S5}
If $\P$, $\Q$ and $\R$ are pre-orders, then, clearly, $\P \times \Q \cong \Q\times \P $ and $(\P \times \Q) \times \R \cong \P \times (\Q \times \R) $
that is, concerning the forcing equivalence of pre-orders, direct product  is a commutative and associative operation. The following lemma generalizes the associativity law.
\begin{lem}   \rm  \label{T4429}
Let $\P$ and $\Q$ be pre-orders and $\la \pi , \leq _\pi , 1_\pi\ra$ a $\P$-name for a pre-order.
Then there is a $\P$-name for a pre-order $\la \pi _1 , \leq _{\pi _1}, 1_{\pi _1}\ra$ such that

(a) $(\P \ast \pi )\times \Q \cong \P \ast \pi _1$;

(b) If $\P$ is $\o$-distributive, $1_\P \Vdash _\P ``\pi $ is $\o _1$-closed" and $\Q$ is $\o _1$-closed, then
$1_\P \Vdash _\P ``\pi _1$ is $\o _1$-closed";

(c) If $1_\P \Vdash _\P ``\pi $ is separative" and $\Q$ is separative, then
$1_\P \Vdash _\P ``\pi _1$ is separative";

(d) If $1_\P \Vdash _\P ``\pi $ is atomless" or $\Q$ is atomless, then
$1_\P \Vdash _\P ``\pi _1$ is atomless".
\end{lem}
\dok
It is easy to show that the triple $\la \pi _1 , \leq _{\pi _1}, 1_{\pi _1}\ra$ works, where

$\pi _1 = \{ \la\la \tau , q \ra \check{\;} , p \ra : p\in \P \land \tau \in \dom \pi \land q\in \Q \land p\Vdash _\P \tau \in \pi \}$,

$\leq _{\pi _1}= \{ \la\la\la \tau _0 ,q_0 \ra , \la \tau _1 ,q_1 \ra \ra \check{\;}, p \ra : p\Vdash _\P \tau _0 , \tau _1 \in \pi \land
\tau _0 \leq _{\pi }\tau _1 \land q_0 \leq _\Q q_1 \}$,

$1_{\pi _1}=\la 1_\pi , 1_\Q \ra\check{\;}$.
\hfill $\Box$
\begin{fac}   \rm  \label{T4430}
Let $\B$ be a non-trivial Boolean algebra, ${\mathcal U}\subset P(\o )$ a non-principal ultrafilter and $\B ^\o /{\mathcal U}$ the corresponding ultrapower. Then

(a) The poset $(\B ^\o /{\mathcal U})^+$ is $\o _1$-closed and separative (folklore);

(b) If the algebra $\B$ is atomless, then $(\B ^\o /{\mathcal U})^+$ is an atomless poset (folklore);

(c) (See \cite{Her}) The poset $(\rp (\B ))^+ $ is forcing equivalent to the two-step iteration $(P(\o )/\Fin )^+ \ast (\B ^\o /\Gamma _1)^+$.
\end{fac}

\begin{te}   \rm  \label{T4421}
For each countable ordinal $\a \geq \o + \o$ the partial order $\la \P (\a ), \subset \ra$ is forcing equivalent to a two-step iteration of the form
$(P(\o )/\Fin )^+ \ast \pi$, where $[\o ] \Vdash$ ``$\pi$ is an $\o _1$-closed, separative atomless forcing".
\end{te}
\dok
Using the notation of Theorem \ref{T4407}, for $\a =\o ^{\g _n +r_n }s_n + \dots + \o ^{ \g _0 +r_0 }s_0 +k $ we have
$\sq \la \P (\a ), \subset \ra \cong \prod _{i=0}^n ( ( \rp ^{r_i}( P(\o ^{\g _i} )/ \I _{\o ^{\g _i} }))^+ )^{s_i}$.

If $r_{i}=0$, for all $i\leq n$, then $\a =\o ^{\g _n  }s_n + \dots + \o ^{ \g _0  }s_0 +k $, where $\g _n \in \Lim$  or $\g _n=1$,
and $\sq \la \P (\a ), \subset \ra \cong \prod _{i=0}^n ( ( P(\o ^{\g _i} )/ \I _{\o ^{\g _i} })^+ )^{s_i}$.
So, if $\g _n\geq \o$, then, by the associativity of direct products,
$\sq \la \P (\a ), \subset \ra \cong  ( P(\o ^{\g _n} )/ \I _{\o ^{\g _n} })^+ \ast \Q$,
where $\Q$ is an $\o _1$-closed, separative and atomless poset (see Theorem \ref{T4121} and Facts \ref{T4401} and \ref{T4043}(a)).
Thus, by Theorem \ref{T4420}, the poset $\sq \la \P (\a ), \subset \ra$ is forcing equivalent to the product
$\R= ((P(\o )/\Fin )^+ \ast \pi) \times \Q$, where $[\o ] \Vdash$ ``$\pi$ is an $\o _1$-closed, separative atomless forcing" and,
by Lemma \ref{T4429}, $\R$ forcing equivalent to an iteration
$(P(\o )/\Fin )^+ \ast \pi _1$, where $[\o ] \Vdash$ ``$\pi _1$ is an $\o _1$-closed, separative atomless forcing".
If $\g _n =1$, then $\a =\o \cdot s_n$ and, by the assumption, $s_n\geq 2$.
Thus $\sq \la \P (\a ), \subset \ra \cong (P(\o )/\Fin )^+ \times ((P(\o )/\Fin )^+)^{s_n -1}=(P(\o )/\Fin )^+ \times \pi$, where
$\pi =(((P(\o )/\Fin )^+)^{s_n -1})\check{\;}$.

If $r_{i_0}>0$, for some $i_0\leq n$, then, by the associativity and commutativity of direct products,
$\sq \la \P (\a ), \subset \ra \cong (\rp ( \rp ^{r_{i_0}-1}( P(\o ^{\g _{i_0}} )/ \I _{\o ^{\g _{i_0}} })))^+ \times \Q$,
where $\Q$ is an $\o _1$-closed, separative and atomless poset (see Theorem \ref{T4121}, Lemma  \ref{T4403} and Fact \ref{T4043}(a)).
If $\B =\rp ^{r_{i_0}-1}( P(\o ^{\g _{i_0}} )/ \I _{\o ^{\g _{i_0}} })$, then , by Fact \ref{T4430}(c), $\sq \la \P (\a ), \subset \ra$ is forcing equivalent
to the product $((P(\o )/\Fin )^+ \ast \pi) \times \Q$, where $[\o ]\Vdash \pi = (\check{\B ^\o } /\Gamma _1)^+$,
and, by Fact \ref{T4430}(a) and (b) applied in extensions by $(P(\o )/\Fin )^+$,
$[\o ] \Vdash$ ``$\pi$ is an $\o _1$-closed, separative atomless forcing". By Lemma \ref{T4429}, $\sq \la \P (\a ), \subset \ra$ is forcing equivalent to an iteration
$(P(\o )/\Fin )^+ \ast \pi _1$, where $[\o ] \Vdash$ ``$\pi _1$ is an $\o _1$-closed, separative atomless forcing".
\hfill $\Box$
\begin{te}   \rm  \label{T4431}
If ${\mathfrak h}=\o _1$, then for each countable ordinal $\a \geq \o $ the partial order $\la \P (\a ), \subset \ra$ is forcing equivalent to
$(P(\o )/\Fin )^+ $.
\end{te}
\dok
If $\a < \o +\o $, then, by Theorem \ref{T4407},  $\sq \la \P (\a ), \subset \ra \cong (P(\o )/\Fin )^+$.

Otherwise, by Theorem \ref{T4421}, $\la \P (\a ), \subset \ra$ is forcing equivalent to a two-step iteration
$(P(\o )/\Fin )^+ \ast \pi$, where $[\o ] \Vdash$ ``$\pi$ is an $\o _1$-closed, separative atomless forcing".
Now, $V\models {\mathfrak h}=\o _1$ implies that CH holds in each generic extension $V_{(P(\o )/\Fin )^+}[G]$ and, by Fact \ref{T4043}(b)
applied in $V_{(P(\o )/\Fin )^+}[G]$, the pre-order
$\pi _G$ is forcing equivalent to $((P(\o )/\Fin )^+)^{V[G]}$. But, since forcing by $(P(\o )/\Fin )^+$ does not produce reals,
$((P(\o )/\Fin )^+)^{V[G]}=((P(\o )/\Fin )^+)^V$ and, hence, $\langle \P (\a ) , \subset \rangle$ is forcing equivalent to
$(P(\o )/\Fin )^+ \times (P(\o )/\Fin )^+$. Now, in $V$ we have ${\mathfrak c}^{<\o _1}={\mathfrak c}$ and the posets $(P(\o )/\Fin )^+$ and
$(P(\o )/\Fin )^+ \times (P(\o )/\Fin )^+ $ are $\o _1$-closed of size ${\mathfrak c}$. In addition, ${\mathfrak h}=\o _1$ implies that
they collapse ${\mathfrak c}$ to $\o _1$ and, by Fact \ref{T4043}(c) they are forcing equivalent (to Coll$(\o _1,{\mathfrak c})$).
\hfill $\Box$
\begin{ex} \rm \label{EX4000}
If $\h _n$ denotes the distributivity number of the poset $((P(\o )/\Fin )^+)^n$, then, clearly,
$\h \geq \h _2 \geq \h _3 \geq \dots \geq \o _1$ and, by Corollary \ref{T4410},
$\h(\sq \la \P (\o n ), \subset \ra ) =\h _n$.
By a result of Shelah and Spinas \cite{SheSpi1}, for each $n\in \N$ there is a model of ZFC in which $\h _{n+1}<\h _n$ and, hence,
the posets $\la \P (\o n ), \subset \ra $ and $\la \P (\o (n+1) ), \subset \ra$ are not forcing equivalent.
\end{ex}
\section{Forcing with quotients over ordinal ideals}\label{S6}
The ideals $\I _{\o ^\d} =\{ I\subset \o ^\d : \o ^\d \not\hookrightarrow I\}$, where $0<\d <\o _1$, are called {\it ordinal}
or {\it indecomposable ideals}. If $\d =\g +r$, where $\g \in \Lim \cup \{ 1 \}$ and $r\in \o$, then, by Facts \ref{T4401}, \ref{T4042} and Theorem \ref{T4407},
we have
\begin{equation}\label{EQ4449}
\sq \la \P (\o ^\d), \subset \ra= (P(\o ^{\g +r})/\I _{\o ^{\g +r}})^+ \cong (\rp ^r(P(\o ^\g)/\I _{\o ^\g }))^+ .
\end{equation}
Let $\h _{\o ^\d} = \h ((P(\o ^{\d})/\I _{\o ^{\d}})^+ )$ and $\t _{\o ^\d} = \t ((P(\o ^{\d})/\I _{\o ^{\d}})^+ )$. Then we have
\begin{te}   \rm  \label{T4439}
For each $\g \in \Lim \cup \{ 1\}$ we have

(a) $\h \geq \h _{\o ^\g}\geq \h _{\o ^{\g +1}}\geq \dots \geq  \h _{\o ^{\g +r}}\geq \dots \geq \o _1$ and, hence, there is $r_0\in \o$ such that
$\h _{\o ^{\g +r}}=\h _{\o ^{\g +r_0}}$, for each $r\geq r_0$;

(b) $\t \geq \t _{\o ^\g}\geq \t _{\o ^{\g +1}}\geq \dots \geq  \t _{\o ^{\g +r}}\geq \dots \geq \o _1$ and, hence, there is $r_0\in \o$ such that
$\t _{\o ^{\g +r}}=\t _{\o ^{\g +r_0}}$, for each $r\geq r_0$.
\end{te}
\dok
(a) By Theorem \ref{T4121}, for each $\d <\o _1$ the poset $\sq \la \P (\o ^\d), \subset \ra$ is $\o _1$-closed and, by Theorem \ref{T4421},
$(P(\o )/\Fin )^+\hookrightarrow _c \sq \la \P (\o ^\d), \subset \ra$. Thus
$\o _1 \leq \t _{\o ^\d}\leq \h _{\o ^\d}\leq \h$.
It is known (see \cite{Kur2}) that $\h ((\rp (\B ))^+)\leq \h (\B ^+)$, for each Boolean algebra $\B$ satisfying $\h (\B ^+)\geq \o _1$, so,
by (\ref{EQ4449}),
$\h _{\o ^{\g +r +1}}
= \h((\rp ^{r+1}(P(\o ^\g)/\I _{\o ^\g }))^+)
=\h((\rp (\rp ^r(P(\o ^\g)/\I _{\o ^\g })))^+)
\leq \h((\rp ^r(P(\o ^\g)/\I _{\o ^\g }))^+)
=\h ((P(\o ^{\g +r})/\I _{\o ^{\g +r}})^+ )$
$=\h _{\o ^{\g +r }}$.

(b) First we prove that $\t _{\o ^\g }\leq \t$, for $\g \in \Lim$. By Proposition \ref{T4425},
$\la P(\o ^\g )\setminus \I _{\o ^\g}, \subset\ra $ $\cong \la P(L)\setminus \I , \subset\ra$
which implies $(P(\o ^\g )/ \I _{\o ^\g})^+\cong ( P(L)/ \I )^+$. Thus, by Claim \ref{T4433},
$\t _{\o ^\g }
=\t ((P(\o ^\g )/ \I _{\o ^\g})^+ )
=\t ( ( P(L)/ \I )^+ )
=\t(P(L)\setminus \I , \subset _\I)
\leq \t$.
The rest of the proof is similar to the proof of (a).
We use the fact (see \cite{Kur2}) that $\t ((\rp (\B ))^+)\leq \t (\B ^+)$, for each Boolean algebra $\B$ satisfying
$\t (\B ^+)\geq \o _1$.
\hfill $\Box$
\begin{ex} \rm \label{EX4001}
By Corollary \ref{T4404}(a) we have $\I _{\o ^2}\cong \Fin \times \Fin$ and, hence,  $\h _{\o ^2}=\h((P(\o \times \o)/(\Fin \times \Fin ))^+)$.
In \cite{Her} Hern\'andez-Hern\'andez proved that in the Mathias model
$\h ((P(\o \times \o )/(\Fin \times \Fin ))^+)=\o _1 $, while $\h ={\mathfrak c} =\o _2$. So, by Theorem \ref{T4439}, in this model
we have $\o _2 ={\mathfrak c}=\h =\h _{\o ^1}>\h _{\o ^2}=\h _{\o ^3}= \dots =\o _1$.

By a result of Szyma\'nski and Zhou \cite{Szym} the poset $(P(\o \times \o)/(\Fin \times \Fin ))^+$
is not $\omega _2$-closed. Thus, by Theorem \ref{T4439}(b), $\t _{\o ^2}=\t _{\o ^3}= \dots =\o _1$ holds in ZFC.
\end{ex}


{\footnotesize
\begin{thebibliography}{99}
%\bibitem{Brend}
%      J.\ Brendle,
%      Distributivity numbers of ${\mathcal P}(\o )/$fin and its friends,
%      Suuri kaiseki kenkyuusho koukyuuroku 1471 (2006) 9--18.
\bibitem{Fra}
      R.\ Fra\"{\i}ss\'{e},
      Theory of relations, Revised edition with an appendix by Norbert Sauer,
      Studies in Logic and the Foundations of Mathematics, 145, North-Holland, Amsterdam, 2000.
\bibitem{Her}
      F.\ Hern\'andez-Hern\'andez,
      Distributivity of quotients of countable products of Boolean algebras,
      Rend.\ Istit.\ Mat.\ Univ.\ Trieste 41 (2009) 27--33 (2010).
\bibitem{Jech}
      T.\ Jech,
      Set Theory, 2nd corr.\ Edition,
      Springer, Berlin, 1997.
\bibitem{Kun}
      K.\ Kunen,
      Set Theory,
      An Introduction to Independence Proofs,
      (North-Holland, Amsterdam, 1980).
\bibitem{KurTod}
      M.\ S.\ Kurili\'c, S.\ Todor\v cevi\'c,
      Forcing by non-scattered sets,
      Ann.\ Pure Appl.\ Logic 163 (2012) 1299--1308.
\bibitem{Kur}
      M.\ S.\ Kurili\'c,
      From $A_1$ to $D_5$: Towards a forcing-related classification of relational structures,
      submitted. http://arxiv.org/abs/1303.2572
\bibitem{Kur1}
      M.\ S.\ Kurili\'c,
      Posets of copies of countable scattered linear orders,
      submitted. http://arxiv.org/abs/1303.2598
\bibitem{Kur2}
      M.\ S.\ Kurili\'c,
      Reduced products,
      submitted.
\bibitem{Lav}
      R.\ Laver,
      On Fra\"{\i}ss\'{e}'s order type conjecture,
      Ann.\ of Math.\ 93,2 (1971) 89-–111.
\bibitem{Rosen}
      J.\ G.\ Rosenstein,
      Linear orderings,
      Pure and Applied Mathematics, 98, Academic Press, Inc.,
      Harcourt Brace Jovanovich Publishers, New York-London, 1982.
\bibitem{SheSpi1}
      S.\ Shelah, O.\ Spinas,
      The distributivity numbers of finite products of $P(\omega)/$fin,
      Fund.\ Math.\ 158,1 (1998) 81--93.
\bibitem{Szym}
      A.\ Szyma\'nski, Zhou Hao Xua,
      The behaviour of $\o ^{2 ^*}$ under some consequences of Martin's axiom,
      General topology and its relations to modern analysis and algebra, V (Prague, 1981), 577–-584,
      Sigma Ser.\ Pure Math., 3, Heldermann, Berlin, 1983.
\end{thebibliography}
}
\end{document}